\documentclass[twoside,bezier,reqno]{amsart}

\usepackage{epsfig}
\usepackage{amsmath,amsthm,amsfonts,amssymb,amscd,euscript,hyperref}

\input{epsf}
\newcommand{\filebegin}{\begin{document}}
\newcommand{\fileend}{\end{document}}
\makeatother
\def\thefootnote{}
\newcommand{\lo}{\longrightarrow}
\newcommand{\NMM}{\hspace*{2mm}}

\renewcommand{\baselinestretch}{1.1}

\input{epsf}
\renewcommand{\baselinestretch}{1.1}
\def\n{\noindent}%

\numberwithin{equation}{section}
\def\mapdown#1{\Big\downarrow\rlap
{$\vcenter{\hbox{$\scriptstyle#1$}}$}}
\newtheorem{theorem}{Theorem}[section]
\newtheorem{lemma}[theorem]{Lemma}
\newtheorem{proposition}[theorem]{Proposition}
\newtheorem{corollary}[theorem]{Corollary}

\newtheorem{conjecture}[theorem]{Conjecture} %---- ADDED COMMAND

\theoremstyle{definition}
\newtheorem{definition}[theorem]{Definition}
\newtheorem{example}[theorem]{\sc Example}
\newtheorem{xca}[theorem]{Exercise}
\theoremstyle{remark}
\newtheorem{remark}[theorem]{Remark}
\begin{document}
%%%%%%%%%%%%%%%%%%%%%%%%%%%%%%%%%%%%%%%

\setcounter{page}{1} \noindent
%Iranian Journal of Mathematical Sciences and Informatics \\
%Vol. x, No. x (201x), pp xx-xx

%%%%%%%%%%%%%%%%%%%%%%%%%%%%%%%%%%%%%%%
\vspace*{2cm}
\begin{center}
{\bf\large On the Closed-Form Solution of a Nonlinear Difference Equation and Another Proof to Sroysang's Conjecture}
 \\[0.5cm]
{Julius Fergy T. Rabago \\[2mm]
Department of Mathematics and Computer Science\\
College of Science\\ 
University of the Philippines Baguio\\ 
Governor Pack Road, Baguio City 2600\\ 
Benguet, PHILIPPINES \\[2mm]
{\tt E-mail: jfrabago@gmail.com, jtrabago@upd.edu.ph}
} \\[2mm]
\end{center}%
\vspace*{0.5cm}
\begin{quotation}
\noindent
{\footnotesize
{\sc Abstract.}
The purpose of this paper is twofold.
First, we derive theoretically, using appropriate transformation on $x_n$, the closed-form solution of the nonlinear difference equation
\[
x_{n+1} = \frac{1}{\pm 1 + x_n},\qquad n\in \mathbb{N}_0.
\]
We mention that the solution form of this equation was already obtained by Tollu et al. in 2013, but through induction principle, and one of our purpose is to clearly explain how was the formula appeared in such structure.
After that, with the solution form of the above equation at hand, we prove a case of Sroysang's conjecture (2013); 
i.e., given a fixed positive integer $k$, we verify the validity of the following claim:
\[
\lim_{x \rightarrow \infty}\left\{ \frac{f(x+k)}{f(x)}\right\}= \phi,
\]
where $\phi=(1+\sqrt{5})/2$ denotes the well-known golden ratio and the real valued function $f$ on $\mathbb{R}$ satisfies the functional equation 
$f(x+2k)=f(x+k) + f(x)$ for every $x\in \mathbb{R}$. 
We complete the proof of the conjecture by giving out an entirely different approach for the other case.}
\end{quotation}

\n \textbf{2000 Mathematics subject classification: } 11B39, 11B37.%\\

\markboth 
{J. F. T. Rabago}
 {Closed-form solution and another proof}

%%%%%%%%%%%%%%%%%%%%%%%%%%%%%%%%%%%%%%%%%%%%%%%%%%%%%%%%%%%%%%%%%%%%%%%%%%%%%%%%%%%%%%%%%%%%%%%%%%%%%%%%%%%%%%%%%%
\newpage
%%%%%%%%%%%% INTRODUCTION %%%%%%%%%%%%%%%
\section{Introduction}

The classical Fibonacci sequence $(F_n)_{n\in\mathbb{N}_0}=\{0,1,1,2,3,5,8,\ldots F_{n+2}=F_{n+1}+F_n,\ldots\}$ has been studied for many years 
and the subject continues to attract attentions of researchers as more fascinating results involving these numbers are obtained and discovered.
In fact, with the growing interest on the topic, various extensions and generalizations of the sequence have been proposed and thoroughly investigated by many mathematicians in the last decades.
Perhaps, the most celebrated of these is the one put forward by A. F. Horadam \cite{horadam} in his 1965 seminal paper, 
from which, we believe, a large amount of recent developments on the topic were based.
In 1878, however, E. Lucas \cite{lucas} first made an extensive investigation of two mere instances of \emph{Horadam sequence} $(w_n)_{n\in\mathbb{N}}=(w_n(w_0,w_1;p,q)):= \{w_0, w_1, \ldots, w_{n+2}=pw_{n+1}-qw_n, \ldots\}$. 
Particularly, Lucas obtained many interesting properties of the sequence $(u_n)_{n\in\mathbb{N}_0}:= u_0=0, u_1=1, u_{n+2}=pu_{n+1}-qu_n$
and $(v_n)_{n\in\mathbb{N}_0}:= v_0=2, v_1=p, v_{n+2}=pv_{n+1}-qv_n$. 
The former is known today as the \emph{fundamental Lucas sequence} and the latter is known as the \emph{primordial Lucas sequence}.
For a good survey of recent developments on Horadam sequences, we refer the readers to \cite{lbf}.

Probably, one of the most interesting property of the Fibonacci sequence is its relation to the widely known golden ratio (or golden number \cite{hambridge, pacioli}) $\phi := (1+\sqrt{5})/2 = 1.6180339887 \ldots$ (cf. Sequence No. A001622 in O.E.I.S.).
It is known (see, e.g., \cite[p. 101]{livio} and \cite[p. 28]{tatterstal}) that the ratios of successive terms of the Fibonacci sequence $(F_{n+1}/F_n)=:(\varphi_n)$  (or any Fibonacci-like sequence) converges to the golden number.
This intriguing value appears vastly in nature and is seen to operate as a universal law in the arrangements of parts such as leaves and branches in plants, and branchings of veins and nerves in animals \cite{padovan1, padovan2}.
This ratio has also been used to analyze the proportions of natural and man-made objects as applications in architecture, designs, paintings, etc \cite{padovan2}. 
Further discussion of this fascinating ratio can be found in \cite{dunlap} and also in \cite{vajda}.
In addition, there is a huge amount of studies scattered in literature about Fibonacci sequence and journals which are entirely devoted to the Fibonacci sequence and its extensions.
For some applications of Fibonacci sequence, we refer the readers to a book by T. Koshy \cite{koshy}.

Recently, numerous papers dealt with various problems relating Fibonacci numbers to other fields of mathematics.
For instance, in \cite{tollu}, Tollu et al. studied the following nonlinear difference equation:
\begin{equation}
	\label{eq1}
	x_{n+1} = \frac{1}{\pm 1 + x_n},\qquad n\in \mathbb{N}_0.
\end{equation}
The authors \cite{tollu} obtained many interesting results regarding the solution of the equation \eqref{eq1} 
and perhaps the most important of these was the closed-form solution of the given equation. 
The solution form of \eqref{eq1}, which we give in the theorem below, was established by Tollu et al. through mathematical induction.

%----- THEOREM 1
\begin{theorem}[\cite{tollu}]
\label{thmtollu}
For any initial value 
$
x_0 \in \mathbb{R}\setminus\left( \left\{ \alpha^{\mp 1},\beta^{\mp 1} \right\} \cup \{\mp \varphi_m\}_{m=1}^{\infty}\right),
$
the closed-form solution of the nonlinear difference equation \eqref{eq1} is given by
\begin{equation}
	\label{sol1}
	x_n = \frac{F_{\pm n}+F_{\pm(n-1)}x_0}{F_{\pm(n+1)}+F_{\pm n}x_0}, \qquad \forall n\in \mathbb{N}_0. 	
\end{equation}
Here $\alpha$ and $\beta$ are the positive and negative root of the quadratic equation $x^2-x-1=0$, respectively (i.e., $\alpha=\phi$ and $\beta=1-\phi$).
\end{theorem}
This result, however, was not supported by any mathematical theory nor was explained by the authors how they were obtained.
Perhaps, this was purposely omitted for some reasons we do not know.

%------------ REMARK 1
\begin{remark}
The exclusion of $\{-\varphi_m\}_{m\in\mathbb{N}}$ and $\{\alpha^{-1},\beta^{-1}\}$ from the set of admissible values for $x_0$, say $A$, 
for the difference equation $x_{n+1}=1/(1+x_n)$ can be explain as follows: 
first, it is evident that $1+x_0$ should not equate to $0$, or equivalently $x_0\neq -1$, so as to have a well-define solution to the equation. 
Now, if for some $n_0 \in\mathbb{N}$, $x_{n_0+1}=-1$, then we must exclude the value of $x_{n_0}$ for which $1/(1+x_{n_0})=-1$ from $A$. 
In this case, $x_{n_0}$ is equal to $-2$ which is again must be removed from $A$.
Continuing the process, we get a sequence of numbers which is in fact the set $\{-\varphi_m\}_{m\in\mathbb{N}}$. 
On the other hand, if $x_0=\alpha^{-1}$ (resp. $x_0=\beta^{-1}$), then we get $x_n=\alpha^{-1}$ (resp. $x_n=\beta^{-1}$) for every $n\in\mathbb{N}$.
This simply means that $\alpha^{-1}$ (resp. $\beta^{-1}$) is a fixed solution of \eqref{eq1}.
A similar explanation holds for the difference equation $x_{n+1}=1/(-1+x_n)$ wherein we must not take any $x_0$ from the set  $\{\varphi_m\}_{m\in\mathbb{N}}$
to have a well-defined solution, and for $x_0=\alpha$ (resp. $x_0=\beta$) we get a fixed solution $x_n=\alpha$ (resp. $x_n=\beta$) for all $n\in\mathbb{N}$ (cf. \cite[Theorem 3]{tollu}).
\end{remark}
%-------------

In an earlier paper, Bacani and the author \cite{bacani} studied a certain generalization of the difference equation \eqref{eq1}.
More precisely, they considered the difference equation 
\begin{equation}
	\label{eq2}
	x_{n+1} = \frac{q}{\pm p +x_n^{\nu}},\qquad n\in\mathbb{N}_0,
\end{equation}
for some initial value 
$
x_0 \in \mathbb{R}\setminus\left( \left\{ \Phi_+^{\mp 1},\Phi_-^{\mp 1} \right\} \cup \{\mp \Psi_m\}_{m\in\mathbb{N}}\right),
$
where $p$ and $q$ are some positive real numbers and $\nu \in\mathbb{N}$.
Moreover, in this case, $\Phi_+$ and $\Phi_-$ represent the positive and negative root of the quadratic equation $x^2-px-q=0$, respectively, and $\Psi_n=w_{n+1}(0,1;p,-q)/w_n(0,1;p,-q)$.
For the case $\nu=1$, the solution form of \eqref{eq2} which is given by
\begin{equation}
	\label{sol2}
	x_n = q\left\{\frac{w_{\pm n}(0,1;p,-q)+w_{\pm(n-1)}(0,1;p,-q)x_0}{w_{\pm(n+1)}(0,1;p,-q)+w_{\pm n}(0,1;p,-q)x_0}\right\}, \qquad \forall n\in \mathbb{N}_0,
\end{equation}
was proved by means of induction principle (cf. \cite[Theorem 1]{bacani}).

Meanwhile, J. S. Han, H. S. Kim and J. Neggers defined in \cite{han} what they called Fibonacci functions 
-- real-valued functions $f$ on $\mathbb{R}$ which satisfy the functional equation $f(x+2)=f(x+1)+f(x)$ for all $x\in\mathbb{R}$.
They developed this notion of Fibonacci function using the concept of $f$-even and $f$-odd functions.
Following \cite{han}, B. Sroysang extended the idea of Fibonacci functions to Fibonacci functions with period $k$
-- real-valued functions on $\mathbb{R}$ which satisfy the equation $f(x+2k)=f(x+k)+f(x)$ for some fixed integer $k\geq 1$, for every $x\in\mathbb{R}$.
Sroysang also defined, in a similar fashion, what he called an odd Fibonacci function with period $k$ 
(a real-valued function which satisfies the equation $f(x+2k)=-f(x+k)+f(x)$, for some fixed integer $k\geq 1$, for all $x\in\mathbb{R}$).
Further, Sroysang made the following conjecture about the asymptotic exponential growth rate of Fibonacci function with period $k$:
\begin{conjecture}[{\cite[Conjecture 25]{sroysang}}]
\label{c1}
If $f$ is a Fibonacci function with period $k\in\mathbb{N}$, then
$
\lim_{x \rightarrow \infty}\left\{ f(x+k)/f(x)\right\}= \phi.
$
\end{conjecture}
These notions of Fibonacci functions were then generalized by the author \cite{rabago} to Horadam functions with period $k$ 
(or second-order linear recurrent functions with period $k$) using the following definition:
\begin{definition}[\cite{rabago}]
\label{def}
Let $k$ be a positive integer, and $r$ and $s$ be positive real numbers. 
A function $w: \mathbb{R} \rightarrow \mathbb{R}$ is said to be a, respectively an odd, second-order linear recurrent function (or simply Horadam function) with period $k$  if it satisfies the functional equation
$w(x+2k)=\pm r w(x+k)+sw(x)$, for every $x\in\mathbb{R}$.
\end{definition}
The above definition naturally gave rise to the concept of Pell and Jacobsthal function (the case $(r,s)=(2,1)$ and $(1,2)$, respectively), 
as well as odd Pell and odd Jacobsthal functions, which are basically analogues of Fibonacci and odd Fibonacci functions, respectively.
Several properties of Horadam functions were studied in \cite{rabago} including the convergence of the ratio $w(x+k)/w(x)$ as $x$ tends to infinity.
This result in fact provides a more general result for the conjecture in \cite{sroysang} and was proven by the author using continued fraction expansion for the root of a non-square integer.
In this note, however, we shall provide a proof of Conjecture \ref{c1} entirely different to the one proposed in \cite{rabago}.
The proof consists of two cases: the first in which we consider the possibility that (i) $f(x+k)/f(x) < 0$, and the other which (ii) $f(x+k)/f(x) \geq 0$.
In the first case we utilize the solution form of \eqref{eq2} and in the second case, we use an entirely different approach which parallels that seen in an elementary analysis course.
As usual, we first prove the existence of the limit of the sequence $(f(x+k)/f(x))$ as $x$ tends to infinity and then show that this limit is nothing but the golden ratio $\phi$. 
%----------------- FLOW OF DISCUSSION

Now we turn on the organization of the rest of the paper.
In the next section (Section 2), we present a theoretical approach in deriving the closed-form solution of the nonlinear difference equation \eqref{eq2}, thus giving a theoretical explanation to Tollu et al.'s result in \cite{tollu}.
In Section 3, we provide another approach in proving Sroysang's Conjecture \ref{c1}. 
The approach we use considers two separate cases. 
In the first case, we utilize the solution form of equation \eqref{eq1} given by \eqref{sol1},
while the second case shall be treated in an entirely different  way.

\section{A Theoretical Approach to Equation \eqref{eq1}}

Consider the nonlinear difference equation given by \eqref{eq2}. 
Obviously, if $\nu = 1$ and $(p,q)=(1,1)$ in \eqref{eq2}, then we'll recover the difference equation \eqref{eq1}.
In this section, we establish the solution form of the difference equation
\begin{equation}
	\label{eq3}
	x_{n+1} =\frac{q}{p+x_n}, \qquad n\in \mathbb{N}_0,
	\tag{T.1}
\end{equation}
through an analytical approach and not with the usual induction method. 
This, in turn, provides a theoretical explanation of the result presented in \cite[Theorem 1]{tollu} concerning the closed-form solution of the given difference equation.
We mention that the same approach can be followed inductively to obtain the solution form of the nonlinear difference equation
\begin{equation*}
	\label{eq4}
	y_{n+1} =\frac{q}{-p+y_n}, \qquad n\in \mathbb{N}_0,
\end{equation*}
so we omit it.

Now we derive the solution form \eqref{sol2} of equation \eqref{eq3} as follows.
We make the substitution $x_n=t_n/t_{n+1}$ (with $n$ replaced by $n-1$) in \eqref{eq3} to obtain the linear  (homogenous) difference equation
\[
\frac{t_n}{t_{n+1}} 
	%= \frac{q}{p+ \dfrac{t_n}{t_{n-1}}}
	= \frac{qt_n}{pt_n+ t_{n-1}}
	\quad \Longleftrightarrow\quad
	t_{n+1}=\frac{p}{q}t_n + \frac{1}{q}t_{n-1}.
\]
In the case when $(t_0,t_1)=(0,1)$, we shall have the solution $t_n=q^{-(n-1)}u_n$, where $u_n$ denotes the $n^{\rm th}$ fundamental Lucas number.
Hence, for arbitrary initial values $(t_0,t_1)$ and after some simple computations, we have 
\[
t_n=t_1 [q^{-(n-1)}u_n] + t_0q^{-1} [q^{-(n-2)}u_{n-1}].
\]
This relation now gives us
\[
x_n = \frac{t_n}{t_{n+1}}
	=\frac{t_1 [q^{-(n-1)}u_n] + t_0 [q^{-(n-1)}u_{n-1}]}{t_1 [q^{-(n-2)}u_n] + t_0[q^{-(n-2)}u_{n-1}]}
	=q\left\{ \frac{u_n + (t_0/t_1) u_{n-1}}{u_{n+1} + (t_0/t_1)u_n}\right\}.
\]
Noting that $x_0=t_0/t_1$ by definition, we get 
\[
x_n=q\left\{ \frac{u_n + x_0 u_{n-1}}{u_{n+1} + x_0u_n}\right\}.
\]
Now since $u_n\equiv w_n(0,1;p,-q)$, we finally have
\begin{equation}
\label{xn}
x_n=q\left\{ \frac{w_n(0,1;p,-q) + w_{n-1}(0,1;p,-q)x_0}{w_{n+1}(0,1;p,-q) + w_n(0,1;p,-q)x_0}\right\},\qquad \forall n\in\mathbb{N}_0,
\tag{T.2}
\end{equation}
which is desired.

%-------- REMARK 2

\begin{remark}
Evidently, for the case when $(p,q)=(1,1)$ in \eqref{eq2}, we recover the result in Theorem \ref{thmtollu}, 
i.e., the closed-form solution of the difference equation $x_{n+1}=1/(1+x_n)$ is given by 
\[
x_n=\frac{w_n(0,1;1,-1) + w_{n-1}(0,1;1,-1)x_0}{w_{n+1}(0,1;1,-1) + w_n(0,1;1,-1)x_0}=\frac{F_n + F_{n-1}x_0}{F_{n+1} + F_nx_0},\qquad \forall n\in\mathbb{N}_0.
\]
\end{remark}
%--------

Now we are in the position to prove the validity of Conjecture \ref{c1} in the next section.

\section{Proof of Sroysang's Conjecture}

%---------- ADDED PART ACCORDING TO REVIEWER
Before we proceed formally with the proof, we first recall the following well-known result in elementary analysis.

\begin{lemma}
\label{lemma}
Let $f$ be a real-valued function continuous on a domain $\mathcal{D} \subset \mathbb{R}$.
Also, let $(x_n)_{n\in \mathbb{N}_0}$ be a convergent sequence in $\mathcal{D}$, with $\lim_{n\to \infty} x_n = \alpha \in \mathcal{D}$.
Then, $\lim_{n\to \infty} f(x_n) = f(\alpha)$.
\end{lemma}
%-----------

%----------- PROOF OF THE CONJECTURE
{\it Proof of the conjecture.} Now, we proceed on proving the conjecture.
Let $k$ be a fixed positive integer and suppose $f$ is a Fibonacci function with period $k$.
Then, $f(x+2k)=f(x+k)+f(x)$ for every $x$ in the real line. 
Note that, for any $x\gg k$, we may write $x$ in the form $\xi+nk$, where $n:=\lfloor x/k\rfloor$.
Denoting $g_n:=\frac{f(\xi+nk)}{f(\xi+(n+1)k)}$, we get
\[
\frac{f(\xi+(n+2)k)}{f(\xi+(n+1)k)}=1+\frac{f(\xi+nk)}{f(\xi+(n+1)k)}
\quad\Longleftrightarrow\quad
g_{n+1}=\frac{1}{1+g_n}.
\]
%-------- CASE 1
\underline{CASE 1}. Suppose first that $g_0 < 0$. 
Note that, in this case, $g_0=f(\xi)/f(\xi+k)$ must not equate to any of the element of the set $\{\beta^{-1}\} \cup \{-\varphi_m\}_{m\in\mathbb{N}}$ 
so as to have a non-fixed and well-defined solution to the nonlinear difference equation $g_{n+1}=1/(1+g_n)$.
Hence, we assume that $g_0\in \mathbb{R}^- \setminus(\{\beta^{-1}\} \cup \{-\varphi_m\}_{m\in\mathbb{N}})$ so that the convergence of the sequence $(g_n)$ may be studied.
Now, using Theorem \ref{thmtollu}, the closed-form solution for $g_n$ is given by
\[
g_n = \frac{F_n+F_{n-1}g_0}{F_{n+1}+F_ng_0}, \qquad \forall n\in\mathbb{N}_0. 
\]
Hence, we can compute for the limit $\lim_{n\rightarrow \infty} \{g_n\}$ as follows:
\begin{align*}
\lim_{n\rightarrow \infty} \{g_n\}
&=\lim_{n\rightarrow \infty} \left\{ \frac{F_n+F_{n-1}g_0}{F_{n+1}+F_ng_0}\right\}
%=\lim_{n\rightarrow \infty} \left\{ \frac{1+\frac{F_{n-1}}{F_n}g_0}{\frac{F_{n+1}}{F_n}+g_0}\right\}\\
=\frac{1+\lim_{n\rightarrow \infty}\left\{\frac{F_{n-1}}{F_n}\right\}g_0}{\lim_{n\rightarrow \infty}\left\{\frac{F_{n+1}}{F_n}\right\}+g_0}\\
&=\frac{1+\phi^{-1}g_0}{\phi+g_0}
=\frac{1}{\phi}.
\end{align*}
Now, in reference to Lemma \ref{lemma}, we get
\begin{align*}
\lim_{n\rightarrow\infty} \left\{ \frac{f(\xi+(n+1)k)}{f(\xi+nk)}\right\} 
&=\lim_{n\rightarrow\infty} \left\{ \frac{1}{g_n}\right\}=\phi.
\end{align*}
However, $x\rightarrow \infty$ as $n\rightarrow \infty$. So we have
\[
\lim_{x\rightarrow\infty} \left\{ \frac{f(x+k)}{f(x)}\right\} 
=\lim_{n\rightarrow\infty} \left\{ \frac{f(\xi+(n+1)k)}{f(\xi+nk)}\right\} 
=\phi.
\]
This proves the first case.\\

%-------- CASE 2
\noindent \underline{CASE 2}. On the other hand, if we assume that $f(x) \geq 0$ and $f(x+k)>0$ for all $x\in\mathbb{R}$ 
(if $f(x) \leq 0$ and $f(x+k) < 0$, then we may define $\hat{f}(x) = -f(x) \geq 0$ and $\hat{f}(x+k)=-f(x+k) > 0$ and then proceed in a similar fashion), 
then we can say that the sequence $(g_n)_{n\in\mathbb{N}_0}$ is well-defined.
We claim that $g_n>0$ for all $n\in\mathbb{N}$. 
To verify this claim, we note that $g_0=f(\xi)/f(\xi+k)\geq 0$ and $g_1=f(\xi+k)/f(\xi+2k)=f(\xi+k)/[f(\xi+k)+f(\xi)]=:M>0$.
Now suppose $g_m>0$ for all $m\leq n \in \mathbb{N}$. Then, $g_{n+1}=1/(1+g_n)>0$. 
By principle of induction, our claim is verified.
Now consider the difference equation
\begin{equation}
\label{gdiff}
g_{n+1}-g_n = -\frac{g_n-g_{n-1}}{(1+g_n)(1+g_{n-1})}.
\end{equation}
Note that the strict inequality $1/(1+\tilde{g}) < 1$ holds for all $\tilde{g} \in (0,\infty)$.
Meanwhile, for any values of $g_1 > 0$, we get $g_2 = 1/(1+g_1) < 1$ which would then implies that $g_n \in (0,1)$.
Hence, we may assume (WLOG) that $g_1 \in (0,1)$ (if not, then we make an adjustment by taking $g_2$ as $g_1$).
Thus, $1-g_1 > 0$ and since $g_n>0$ for all $n\in \mathbb{N}$, we get
\begin{align*}
	1-g_1+ g_{n-1} > 0\ (n\geq 2) &\qquad \Longleftrightarrow\qquad 1+ g_{n-1} > g_1\\
		&\qquad \Longleftrightarrow\qquad 1+ g_{n-1} + g_n(1+g_{n-1}) > 1 + g_1\\
		&\qquad \Longleftrightarrow\qquad (1 + g_n) (1+ g_{n-1} ) > 1 + g_1.
\end{align*}
Thus (WLOG) we may claim that, for all $n\in\mathbb{N}$, we have the strict inequality 
\[
(1+g_n)(1+g_{n-1})=1+g_{n-1}+g_n+g_{n-1}g_n
> 1+ M,
\]
for all $n\in\mathbb{N}\setminus\{1\}$. 
Taking the absolute value on both sides of equation \eqref{gdiff} now gives us the relation
\[
|g_{n+1}-g_n| 
= \frac{|g_n-g_{n-1}|}{|(1+g_n)(1+g_{n-1})|}
<\frac{|g_n-g_{n-1}|}{1+M}
\]
for all $n\in\mathbb{N}\setminus\{1\}$.
Now, since 
\begin{align*}
|g_2-g_1|
&=\left| \frac{f(\xi+2k)}{f(\xi+3k)}-\frac{f(\xi+k)}{f(\xi+2k)}\right|
=\left| \frac{f(\xi+k)+f(\xi)}{2f(\xi+k)+f(\xi)}-\frac{f(\xi+k)}{f(\xi+k)+f(\xi)}\right|\\
&=\frac{|f^2(\xi)+f(\xi)f(\xi+k)-f^2(\xi+k)|}{(2f(\xi+k)+f(\xi))(f(\xi+k)+f(\xi))}=:c,
\end{align*}
the first iteration gives us 
$
|g_3-g_2| < |g_2-g_1|(1+M)^{-1} =: c(1+M)^{-1}
$
which in turn leads us to
$
|g_4-g_3| < |g_3-g_2|(1+M)^{-1} < c(1+M)^{-2}
$.
Continuing the process up to some integer $n \in\mathbb{N}\setminus \{1\}$, we obtain
$|g_{n+1}-g_n| < c(1+M)^{-(n-1)}$,
which can be verified easily by induction.
Indeed, given the assumption that $|g_{m+1}-g_m|<c(1+M)^{-(m-1)}$, for all $m \leq n$, we have
$|x_{n+2}-x_{n+1}|<(1+M)^{-1}|x_{n+1}-x_n|<c(1+M)^{-n}$.
Thus, $|g_{n+1}-g_n| < c(1+M)^{-(n-1)}$ for all $n \in\mathbb{N}\setminus \{1\}$.

Next, we show the existence of the limit of the sequence $(g_n)_{n\in\mathbb{N}_0}$; that is, we prove that $(g_n)_{n\in\mathbb{N}_0}$ is Cauchy.
To do this, we first approximate the value $|g_m-g_n|$ for arbitrary choice of index $m$ and $n$ (with $m>n$) 
and then show that, for some sufficiently large $N$, $|g_m-g_n|< \varepsilon$ for each $m>n\geq N $, for every $\varepsilon >0$.

We express $g_m-g_n$ as $g_m-g_n=(g_m-g_{m-1})+(g_{m-1}-g_{m-2})+\ldots+(g_{n+2}-g_{n+1})+(g_{n+1}-g_{n})$.
Hence, by triangle inequality, $|g_m-g_n|\leq |g_m-g_{m-1}|+|g_{m-1}-g_{m-2}|+\ldots+|g_{n+2}-g_{n+1}|+|g_{n+1}-g_{n}|$.
Since $|g_{n+1}-g_n|<c(1+M)^{-(n-1)}$, then using the formula for the sum of a geometric series, we now have
\begin{align*}
|g_m-g_n|
&\leq \frac{c}{(1+M)^{m-2}} + \frac{c}{(1+M)^{m-3}}+\ldots+\frac{c}{(1+M)^{n-2}}+\frac{c}{(1+M)^{n-1}}\\
&<\dfrac{\frac{c}{(1+M)^{n-1}}}{1-\frac{1}{1+M}}=\frac{c}{(1+M)^{n-2}}=:\Omega(n).
\end{align*}
Now, given $\varepsilon>0$, we choose a sufficiently large $N$ such that $\Omega(n) < \varepsilon$.
So, for all $m>n\geq N$, $|x_m-x_n|<\Omega(n)\leq\Omega(N)<\varepsilon$.
This proves that $(g_n)_{n\in\mathbb{N}_0}$ is Cauchy, thereby implying that $L:=\lim_{n\rightarrow\infty} g_n$ exists.
Going back to the relation $g_{n+1}=1/(1+g_n)$, we have
\[
L=\lim_{n\rightarrow\infty} \{g_{n+1}\} 
=\lim_{n\rightarrow\infty} \left\{ \frac{1}{1+g_n}\right\}
=\frac{1}{1+ \lim_{n\rightarrow\infty} g_n}
=\frac{1}{1+L}.
\]
This yields the quadratic equation $L^2+L-1=0$ whose solution set is $\{-\phi,\phi-1\}$.
However, we have shown that $g_n$ is positive for every $n\in\mathbb{N}$, 
so $L=\phi-1$. By virtue of Lemma \ref{lemma}, it follows that 
\begin{align*}
\lim_{n\rightarrow\infty} \left\{ \frac{f(\xi+(n+1)k)}{f(\xi+nk)}\right\} 
&=\lim_{n\rightarrow\infty} \left\{ \frac{1}{g_n}\right\}=\frac{1}{L}=\frac{1}{\phi-1}=\phi.
\end{align*}
But $x\rightarrow \infty$ as $n\rightarrow \infty$. Hence, the above equation is equivalent to
\[
\lim_{x\rightarrow\infty} \left\{ \frac{f(x+k)}{f(x)}\right\} 
=\lim_{n\rightarrow\infty} \left\{ \frac{f(\xi+(n+1)k)}{f(\xi+nk)}\right\} 
=\phi,
\]
proving the second case. This completes the proof of the conjecture.\\
%---- END OF PROOF

Now we state Conjecture \ref{c1} as a theorem.
\begin{theorem}
Let $k$ be a positive integer. If $f:\mathbb{R}\rightarrow \mathbb{R}$ is a Fibonacci function with period $k$, then
$
\lim_{x \rightarrow \infty}\{f(x+k)/f(x)\}= \phi.
$
\end{theorem}

%-------- REMARK 3
\begin{remark}
We mention that the above theorem can be proven using Theorem \ref{thmtollu}, irrespective of the sign of the initial value $g_0=f(\xi+k)/f(\xi)$. 
That is, as long as we are sure that $g_0 \in \mathbb{R}\setminus(\{\beta^{-1}\} \cup \{\varphi_m\}_{m\in\mathbb{N}})$, then we know that $g_n$ converges to $\phi^{-1}$ (cf. \cite[Theorem 4]{tollu}).
\end{remark}
%--------

%-------- REMARK 4
\begin{remark}
\label{remark4}
Obviously, the same approach can be applied to prove a more general result of Conjecture \ref{c1}.
More specifically, we can prove that the ratio of Horadam functions $w(x+k)/w(x)$ (a real-valued function satisfying the functional equation $w(x+2k) = rw(x+k) + sw(x)$ for some positive real numbers $r$ and $s$) will converge to the positive root of the quadratic equation $x^2-rx-s=0$ (cf. \cite[Corollary 6.3]{rabago}) 
using the closed-form solution \eqref{sol2} of the case $\nu=1$ of the nonlinear difference equation \eqref{eq2} (cf. \cite[Theorem 1]{bacani}).
That is, the limit $\lim_{x\rightarrow \infty} \{w(x+k)/w(x)\}=(r+\sqrt{r^2+4s})/2=:\rho$ (the positive root of the equation $x^2-rx-s=0$) can be verified as follows:

Denote $h_n:=w(\xi+nk)/w(\xi+(n+1)k)$ so that the functional equation given by $w(x+2k)=rw(x+k)+sw(x)$ (which $w$ satisfies) is transformed into the nonlinear difference equation $h_{n+1}=1/(r+s\cdot h_n)$, for all $n\in\mathbb{N}_0$.
Using the substitution $p=r/s$ and $q=1/s$, we then obtain the nonlinear difference equation $h_{n+1}=q/(p+h_n)$.
Assume that $h_0 \in \mathbb{R}\setminus(\{\Phi_-^{-1}\}\cup\{-u_{m+1}/u_m\}_{m\in\mathbb{N}})$ (where $u_n$ denotes the $n^{th}$ fundamental Lucas numbers).
Then, in view of equation \eqref{xn}, we obtain the closed-form solution
\[
h_n = q\left\{\frac{u_n+u_{n-1}h_0}{u_{n+1}+u_nh_0} \right\}, \qquad \forall n\in\mathbb{N}_0.
\]  
Letting $n\rightarrow \infty$, we get
\begin{align*}
\lim_{n\rightarrow \infty}\{h_n\} 
&= q\left\{\frac{1+\lim_{n\rightarrow \infty}\left\{\dfrac{u_{n-1}}{u_n}\right\}h_0}{\lim_{n\rightarrow \infty}\left\{\dfrac{u_{n+1}}{u_n}\right\}+h_0} \right\}
= q\left\{\frac{1+\Phi_+^{-1}h_0}{\Phi_++h_0}\right\}
=\frac{q}{\Phi_+}.
\end{align*} 
Using the fact that $x\rightarrow \infty$ as $n\rightarrow \infty$, together with Lemma \ref{lemma}, we get
\[
\lim_{x\rightarrow\infty} \left\{ \frac{w(x+k)}{w(x)}\right\} 
=\lim_{n\rightarrow\infty} \left\{ \frac{w(\xi+(n+1)k)}{w(\xi+nk)}\right\} 
=\lim_{n\rightarrow\infty} \left\{ \frac{1}{h_n}\right\} 
=\frac{\Phi_+}{q}.
\]
Note, however, that
\[
\frac{\Phi_+}{q}
=\frac{1}{2q}(p+\sqrt{p^2+4q})
=\frac{s}{2}\left(\frac{r}{s}+\sqrt{\frac{r^2}{s^2}+\frac{4}{s}}\right)
=\rho.
\]
Thus, $\lim_{x\rightarrow\infty} \{ w(x+k)/w(x)\} =\rho$.
\end{remark}
%--------

%-------- REMARK 5
\begin{remark}
We also emphasize that the method used previously to prove that $\lim_{x\rightarrow\infty} \{w(x+k)/w(x)\}=\rho$ 
can definitely be applied to show that the ratio $\varpi(x+k)/\varpi(x)$ of odd Horadam functions with period $k$ 
(satisfying the functional equation $\varpi(x+2k)=-\varpi(x+k)+\varpi(x)$) will converge to $-\rho$ (cf. \cite[Corollary 6.7]{rabago}). 
In this case, the closed-form solution of the nonlinear difference equation $y_{n+1}=q/(-p+y_n)$ (with initial value $y_0 \in\mathbb{R}\setminus(\{\Phi_+\}\cup \{u_{m+1}/u_m\}_{m\in\mathbb{N}}) $) given by
\[
y_n=q\left\{ \frac{w_{-n}(0,1;p,-q) + w_{-(n-1)}(0,1;p,-q)y_0}{w_{-(n+1)}(0,1;p,-q) + w_{-n}(0,1;p,-q)y_0}\right\}, \qquad\forall n\in\mathbb{N}_0
\]
can be utilized. This in turn will prove, as a special case (the instance $(r,s)=(1,1)$), Sroysang's second conjecture:
if $f$ is an odd Fibonacci function with period $k\in\mathbb{N}$, then $\lim_{x \rightarrow \infty}\{f(x+k)/f(x)\}= -\phi$ (cf. \cite[Conjecture 26]{sroysang}). 
 \end{remark}
%--------

%-------- REMARK 6
\begin{remark}
As for our final remark, we mention that the following statement is also true:
\[
\lim_{x\rightarrow -\infty} \left\{ \frac{w(x+k)}{w(x)}\right\} 
	= -\lim_{x\rightarrow -\infty} \left\{ \frac{\varpi(x+k)}{\varpi(x)}\right\}
	=-\rho,
	%\quad \text{and}\quad
	\]
%and
%\[	
%\lim_{x\rightarrow -\infty} \left\{ \frac{\varpi(x+k)}{\varpi(x)}\right\} 
%	=\rho,
% \]
where $w$ and $\varpi$ are Horadam and odd Horadam functions with period $k$, respectively (cf. \cite{rabago}).
Particularly, if $f$ is a Fibonacci function (resp. an odd Fibonacci function) with period $k$, then the sequence of ratios 
$\left\{f(x+k)/f(x)\right\}$ converges to $-\phi$ (resp. $\phi$) as $x$ decreases without bound.
These results can be verified easily with the same approach as above and using the fact that the Horadam numbers, in general, can naturally be extended to negative numbers using the relation 
$w_{-n}(w_0,w_1;p,q)=(-1)^{n+1}w_n(w_0,w_1;p,q)$ together with the solution form of the nonlinear difference equation $h_{n+1}=q/(-p+h_n)$.
\end{remark}
%--------

{\bf Author's Note} It was pointed out by one of the referee of this paper that the difference equation 
\[
	g_{n+1} = \frac{1}{1+g_n} \quad \text{for every $n\in \mathbb{N}_0$}
\]
is related to continued fractions. 
Indeed, for sufficiently large $N > 0$, we may iterate the right hand side of the above equation to obtain
\[
g_N = \dfrac{1}{1+\dfrac{1}{1+\dfrac{1}{1+\dfrac{1}{1+\dfrac{1}{\quad \cdots + \dfrac{1}{1+\frac{1}{g_0}}}}}}}.
\]
Recall that (see, e.g., \cite{rabago1}) $\phi - 1 = [0;1,1,1,\ldots]$ where
\[
 [0;1,1,1,\ldots] = \dfrac{1}{1+\dfrac{1}{1+\dfrac{1}{1+\dfrac{1}{1+\cdots}}}}.
\]
Thus, for every $\varepsilon > 0$, we can find an integer $N>0$, sufficiently large, such that $|g_n - (\phi - 1)|< \varepsilon$.
Equivalently, we have
$
\lim_{N\to \infty} g_N = \phi - 1.
$
Another important thing to note regarding the sequence $(g_n)_{n \in \mathbb{N}}$ (with $g_1$ in the unit interval $(0,1)$) is that, 
the $n$-th term $g_n$ is either contained in the interval $\left[ F_n/F_{n+1}, F_{n+1}/F_{n+2}\right]$ or in $\left[ F_{n+1}/F_{n+2}, F_n/F_{n+1}\right]$ (depending on the parity of $n$) (cf. \cite[Lemma 2.1]{hakami} and \cite[Lemma 5]{rabago1}).
Noting that $F_{n+1}/F_n \to \phi$ as $n \to \infty$, 
one can immediately see (possibly through Cantor's Intersection Theorem \cite{cantor}) that $g_n \to 1/\phi = \phi -1$.
It is worth mentioning that this approach was in fact used explicitly by the author \cite{rabago} to prove a more general case of Sroysang's conjecture (cf. Remark \ref{remark4} above).
%-------- SUMMARY
\section{Summary}

We have verified affirmatively, in an alternative fashion, Sroysang's conjecture regarding the asymptotic growth rate of the so-called Fibonacci functions (and odd Fibonacci functions) with period $k$. 
The technique we have used in proving the conjecture, which is one of the main objective of our work, utilizes some well-known results and direct computations, using elementary properties of classical analysis.
In the proof, we have started with the transformation $f(x)/f(x+k)=f(\xi+nk)/f(\xi+(n+1)k)=:g_n$ with $f$ satisfying the functional equation $f(x+2k)=f(x+k)+f(x)$ for all $x\in\mathbb{R}$,
then utilized the closed-form solution of the difference equation $g_{n+1}=1/(1+g_n)$. 
In this approach, we first showed that the sequence $(g_n)$ converges to $\phi -1$ and then used this fact to arrive at the conclusion that the asymptotic exponential growth rate of Fibonacci function with period $k$ indeed converges to the well-known golden ratio $\phi$.
As a remark, we have also asserted that the same approach can be followed inductively to prove a more general case of the statement. 
Further, it was noted that the idea behind the method used to establish the main result can be employed to verify a similar result for odd Horadam functions with period $k$. 
The desired result for this case, as we have remarked, can be achieved using a property of Horadam numbers with negative indices combined with the solution form of the nonlinear difference equation $h_{n+1}=q/(-p+h_n)$. 
The resulting property, in turn, validates (as a special case) Sroysang's second conjecture.

\section*{Acknowledgments}
The author wishes to thank the anonymous referees for carefully handling and examining his manuscript.
Their constructive comments and suggestions greatly improved the quality of the paper.
The proof of the main result was substantially refined due to the valuable suggestion of one of the referee.% who pointed out Lemma \ref{lemma}.
\newpage	%-----------------enter new page
%%%%%%%%%%%%%%%%%%%%%%%%%%%%%%%%%%%%%%%%%%%%%%%%%%%%%%%%%%
\providecommand{\bysame}{\leavevmode\hbox
to3em{\hrulefill}\thinspace}

%%%%%%%%%%%%%%%%%%%%%%%%%%%%%%%%%%%%%%%%%%%%%%%%%%%%%%%%%%%%%%%%%%%%%%%%%


\begin{thebibliography}{99}

\bibitem{bacani}
J. B. Bacani and J. F. T. Rabago, 
{\it On two nonlinear difference equations},
submitted.

\bibitem{dunlap}
R. A. Dunlap,
{\it The Golden Ratio and Fibonacci Numbers},
World Scientific, Singapore, 1998.

\bibitem{hakami}
A. Hakami, 
{\it An application of Fibonacci sequence on continued fractions},
Int. Math. Forum, {\bf 10}(2), 69--74.

\bibitem{hambridge}
J. Hambidge, 
{\it Dynamic Symmetry: The Greek Vase}, 
New Haven CT: Yale University Press, 1920.

\bibitem{han}
J. S. Han, H. S. Kim, J. Neggers,
{\it On Fibonacci functions with Fibonacci numbers}, 
Adv. Differ. Equ., {\bf 2012} (2012), Article 126, 7 pages. 

\bibitem{horadam}
A. F. Horadam, 
{\it Basic properties of a certain generalized sequence of numbers}, 
Fib. Quart., {\bf 3} (1965), 161--176.

\bibitem{koshy}
T. Koshy,
{\it Fibonacci and Lucas Numbers with Applications}, 
Pure and Applied Mathematics, Wiley-Interscience, New York, 2001.

\bibitem{lbf}
P. J. Larcombe, O. D. Bagdasar, and E. J. Fennessey, 
{\it Horadam sequences: a survey},
Bulletin of the I.C.A., {\bf 67 } (2013), 49--72.

\bibitem{livio}
M. Livio,
{\it The Golden Ratio: The Story of Phi, the World's Most Astonishing Number}, 
New York: Broadway Books, 2002.

\bibitem{lucas}
E. Lucas, {\it Th\'{e}orie des Fonctions Num\'{e}riques Simplement P\'{e}riodiques}, 
American Journal of Mathematics, {\bf 1} (1878), 184--240, 289--321; 
reprinted as ``The Theory of Simply Periodic Numerical Functions'', Santa Clara, CA: The Fibonacci Association, 1969.

\bibitem{oeis}
O.E.I.S. Foundation Inc. (2011), The On-Line Encyclopedia of Integer Sequences, \url{http://oeis.org.}

\bibitem{pacioli}
L. Pacioli, Luca, 
{\it De divina proportione}, 
Luca Paganinem de Paganinus de Brescia (Antonio Capella) 1509, Venice.

\bibitem{padovan1} 
R. Padovan, 
{\it Proportion}, 
Taylor \& Francis, pp. 305--306, 1999. 

\bibitem{padovan2}
R. Padovan, 
{\it Proportion: Science, Philosophy, Architecture}, 
Nexus Network Journal {\bf 4}(1) (2002), 113--122.

%\bibitem{noe}
%Noe, Tony; Piezas, Tito III; and Weisstein, Eric W. 
%``Fibonacci n-Step Number.'' From{ {\it MathWorld}--A Wolfram Web Resource.} 
%\texttt{http://mathworld.wolfram.com/Fibonaccin-StepNumber.html}

\bibitem{rabago}
J. F. T. Rabago, 
{\it On second-order linear recurrent functions with period $k$ and proofs to two conjectures of Sroysang},
Hacet. J. Math. Stat., {\bf 45}(2) (2016), 429--446.

\bibitem{rabago1}
J. F. T. Rabago, 
{\it On k-Fibonacci Numbers with Applications to Continued Fractions},
Journal of Physics: Conference Series {\bf 693} (2016) 012005.

\bibitem{cantor}
T. Rowland, {\it Cantor's Intersection Theorem}. 
From MathWorld--A Wolfram Web Resource, created by Eric W. Weisstein. 
\url{http://mathworld.wolfram.com/CantorsIntersectionTheorem.html}

\bibitem{sroysang}
B. Sroysang,
{\it On Fibonacci functions with period $k$},
Discrete Dyn. Nat. Soc., {\bf 2013} (2013), Article ID 418123, 4 pages.

\bibitem{tatterstal}
J. J. Tattersall,
{\it Elementary Number Theory in Nine Chapters (2nd ed.),} 
Cambridge University Press, 2005. 

\bibitem{tollu}
D. T. Tollu, Y. Yazlik and N. Taskara,
{\it On the solutions of two special types of Riccati difference equation via Fibonacci numbers}, 
Adv. Differ. Equ., {\bf 2013} (2013), Article 174, 7 pages.

\bibitem{vajda}
S. A. Vajda,
{\it Fibonacci \& Lucas Numbers and The Golden Section: Theory and Applications},
Ellis Horwood Ltd., Chishester, 1989.

\end{thebibliography}
\end{document}